\documentclass[12pt]{article}
\usepackage{graphicx,color}
\usepackage{amsfonts}
\usepackage{amssymb}

 \advance \textheight by \topskip
\newtheorem{theorem}{Theorem}
\newtheorem{lemma}{Lemma}

\newcommand{\fx}{{f(x)}}
\renewcommand{\phi}{{\varphi}}
\newcommand{\R}{\mbox {${\rm {I\!R}}$}}

\begin{document}
\title{Limit Sets of Convex non Elastic Billiards}
\author{Roberto Markarian $^1$, Sylvie Oliffson Kamphorst$^2$ and S\^onia Pinto-de-Carvalho$^2$}

\date{}

\maketitle

\begin{abstract}
Inspired by the work of Pujals and Sambarino on dominated splitting \cite{PS09}, we present billiards with a modified reflection law which constitute simple examples of dynamical systems with limit sets having dominated splitting and where the dynamics is a rational or irrational rotation.
\end{abstract}

%\ams{37D30, 37C70}

\section{Introduction}

Let  $M$ be a Riemannian manifold and  $f: M \to M^{\prime}\subset M$ a
diffeomorphism. An $f$-invariant set $\Lambda$  has a dominated splitting if
if its tangent bundle $T_\Lambda M=U\oplus V$, where $U$ and $V$ are non trivial invariant continuous subbundles
such that, for constants  $C>0$ and $0<\gamma<1$:
\begin{equation}
\|Df^{n}{|_{U(x)}}\|\ \|Df^{-n}{|_{V(f^{n}({x}))}}\|^{}\le C \gamma^n,
\mbox{ for
all } x\in\Lambda, n\ge 0.  \label{fds}
\end{equation}
Clearly, any hyperbolic splitting is a dominated one.

Some important consequences of this property on the dynamics were given by Pujals and Sambarino in \cite{PS09}. A spectral decomposition theorem was obtained for $C^{2}$ compact surface
diffeomorphisms having dominated splitting over the limit set $\displaystyle L(f) = \overline{\bigcup _{x\in M}\left(\omega(x) \cup \alpha(x)\right)}$
where $\omega(x)$ and $\alpha(x)$ are the $\omega$ and $\alpha$-limit sets of $x$, respectively:

\textbf{Theorem (Pujals-Sambarino \cite{PS09}):} \emph{Let $M$ be a compact 2-manifold and $f:M \to M^{\prime}\subset M$ a $C^2$-diffeomorphism. Assume that $L(f)$ has a
dominated splitting. Then $L(f)$ can be decomposed into $L(f)=\mathcal{I}\cup{\mathcal{R}}\cup
{\tilde{\mathcal{L}}}$ such that:}

\begin{enumerate}
\item \emph{$\mathcal{I}$ is a set of periodic points with bounded periods
contained in a disjoint union of finitely many normally hyperbolic periodic
arcs or simple closed curves. }
\item \emph{$\mathcal{R}$ is a finite union of normally hyperbolic periodic
simple closed curves supporting an irrational rotation. }
\item \emph{$\tilde{\mathcal{L}}$ can be decomposed into a disjoint union
of finitely many compact invariant and transitive sets (called basic sets).
The periodic points are dense in $\tilde{\mathcal{L}}$ and at most
finitely many of them are non-hyperbolic periodic points. The (basic) sets
above are the union of finitely many (nontrivial) homoclinic classes.
Furthermore $f|_{\tilde{\mathcal{L}}}$ is expansive. }
\end{enumerate}

Our purpose, in this work, is to construct simple examples of dynamical systems with attractors admitting a dominated splitting and where the dynamics is of type $\mathcal{I}$ or $\mathcal{R}$. We will follow the ideas developed in \cite{MPS09} and \cite{AMS09}. In  \cite{MPS09}, non conservative billiards were introduced by a modification of the reflection rule and  the existence of attractors was demonstrated for a wide class of dispersing and semi-dispersing billiards and of billiards with focusing components. In \cite{AMS09} models were studied numerically and different attractors, periodic and chaotic are presented.

We concentrate on  strictly convex billiard tables, with boundary formed by a unique sufficiently differentiable focusing component. They present structures like KAM stability islands and invariant rotational curves non homotopic to a point. We will investigate how, in the presence of non conservative perturbations, an invariant curve will give rise to an attractor.
The maps we consider here are more general than the pinball billiards introduced in \cite{MPS09}, as the perturbation of the angle is not necessarily biased to the normal direction.

The contents of the present paper is outlined in four sections that follow this introduction.
{In} section~\ref{sec:cones} we present the main tools needed to work with dominated splitting while section~\ref{sec:classic} deals with the basic properties of classical billiards on ovals.

In section~\ref{sec:nonela} we introduce our {\sl non elastic billiards}. They are defined as a composition of a classical billiard followed by a change of the reflection angle, corresponding to a contraction in the vertical fibers of an invariant rotational curve.
We will prove that under some differentiability hypotheses and some bounds on the contraction, there exists a compact strip in the phase space,  such that the non elastic billiard map is a $C^2-$diffeomorphism from that strip onto its image. Its limit set contains the invariant curve and has a dominated splitting. Moreover, the non elastic dynamics on the invariant curve is determined by its rotation number with respect to the original classical billiard map.

This result will guide us to construct our examples of non elastic billiards on ovals with dominated splitting and attractors supporting a rational or an irrational rotation. They are presented in section \ref{sec:exa}, where we explore their properties theoretically and numerically.

Our result (Theorem \ref{nonelastic}) could clearly be established in the more general setting of $C^2$ conservative twist maps. Once we have  a $C^2$ invariant rotational curve and assuming the necessary bounds, everything would work likewise in the proof. The main problem is to build specific examples, other than the obvious twist integrable case (which corresponds to the circular billiard) or the oval billiards presented here, and check, for instance and at least numerically, the size of the basin of attraction.

\section{Dominated Splitting, Cone Fields and Quadratic Forms}
\label{sec:cones}

The existence of a dominated splitting follows from the existence of an eventually strictly invariant cone field and of an uniform control of expansions and/or contractions as showed in \cite{Wo01}. The criteria presented there in Proposition~4.1 can be translated, in dimension 2, into the operational lemma below.

Let  $u,v:M\mapsto  TM$ be two vector fields such that  $u(x)=u_x$ and $v(x)=v_x$ are linearly independent vectors in $ T_xM$.
They induce a nondegenerate quadratic form $Q$ on $ TM$ by $Q_x (a u_x + b v_x) = ab$ and a cone field given at each $x$ by $\displaystyle  {\cal C}(x) = \left\{ w \in {T}_xM  \, : \,  Q_x (w) > 0 \right\} \cup \{ 0 \}$ and whose boundaries, at each point, are given by ${\cal C}_{0}(x) = \left\{ w \in {T}_xM  \, : \,  Q_x (w) = 0 \right\}$.
If the vector fields are continuous, the quadratic form and the cone field are also continuous.

Given $x \in M$ and a vector $w = a u_x + b v_x \in {T}_xM$, let
$Df_x w = a_1 u_{\fx} + b_1 v_{\fx} \in {T}_{\fx}M$ denote the image of $w$ under the
derivative $Df_x$.  Then we have
\begin{equation}
\left(
\begin{array}{c}
a_1  \\ b_1
\end{array}
\right)
=
[Df_x]_U
\left(
\begin{array}{c}
a  \\ b
\end{array}
\right)
\label{eq:ab}
\end{equation}
where $[Df_x]_U$ is the matrix representation of the derivative at $x$, with the choice of $\{u_x,v_x\}$ and $\{u_{\fx},v_{\fx}\}$ as bases of ${T}_xM$ and  ${T}_{\fx}M$ respectively.

\begin{lemma}\label{lema1}
Let $\Lambda$ be a compact $f$-invariant subset of $M$. If
there is a choice of vector fields $u,v$ such that the entries of $[Df_x]_U$ are strictly positive for every $x\in \Lambda$, then $\Lambda$ has a dominated splitting.
\end{lemma}
Proof: If the entries of $[Df_x]_U$ are strictly positive for every $x\in \Lambda$ then for every $w = a u_x + b v_x $, $x\in \Lambda$, $ a b \ge 0$, $a^2 + b^2 >0$ we have $a_1 b_1 > 0$ where $Df_x w = a_1 u_{\fx} + b_1 v_{\fx}$.
This implies that for every $x\in \Lambda$ $Df_x ({\cal C}(x) \cup {\cal C}_{0}) \subset {\cal C} (\fx)$, i.e., $\cal C$ is strictly $Df$-invariant ($f$ is strictly $Q$-separated).
It follows from Proposition 4.1 in \cite{Wo01} that $\Lambda$ has a dominated splitting. \hfill  $\blacksquare$

\section{Classical Billiards on Ovals}\label{sec:classic}

Let $\Gamma$ be an oval, i.e., a  plane,  simple, closed, $C^k, k
\ge 3$, curve, with strictly positive curvature, parameterized  counterclockwise by
$\phi$, the angle between the tangent vector and an horizontal
axis. Let $R(\varphi)$ be its radius of curvature at $\varphi$.

The classical billiard problem on $\Gamma$ consists of the free
motion of a point particle inside  $\Gamma$, making elastic reflections at the impacts with the boundary. The
motion is then determined  by the point of reflection at $\Gamma$ and the direction of motion immediately after each reflection.
They can be given by the parameter $\varphi\in [0,2\pi)$, that
will locate the point of reflection and by the angle $\alpha \in (0,\pi)$ between the tangent vector and the outgoing trajectory,  measured counterclockwise.

The classical billiard defines a map $B$ from the open cylinder $[0,2\pi)\times(0,\pi)$ into itself, which has some very well known properties (see, for instance, \cite{kat} and \cite{kato} for the properties of billiards and twist maps listed in this section).

Denoting $B(\varphi_0,\alpha_0)= (\varphi_1,\alpha_1)$, the derivative of $B$ at $(\varphi_0,\alpha_0)$ is
\begin{equation}\label{eq:deriv}
 DB_{(\phi_0,\alpha_0)} =
\frac{1}{R_1 \sin \alpha_1}
\left(
\begin{array}{cc}
L -  R_0 \sin \alpha_0  & L   \\
L -  R_0 \sin \alpha_0 -  R_1 \sin \alpha_1 & L -  R_1 \sin \alpha_1
\end{array}
\right)
\end{equation}
where $R_i = R(\phi_i)$  and $L$ is the distance between $\Gamma(\varphi_0)$ and $\Gamma(\varphi_1)$.

If $\Gamma $ is $C^k$, then $B$ is a $C^{k-1}$-diffeomorphism. It preserves the measure $d\nu=R(\varphi)\sin\alpha\, d\alpha d\varphi$. 

It is reversible with respect to the reversing symmetry $S(\varphi,\alpha)=(\varphi,\pi-\alpha)$, as $S\circ B= B^{-1}\circ S$, meaning that the phase space is symmetric with respect to the line $\alpha = \pi/2$.

Moreover, as $\Gamma$ is an oval, $B$  has the monotone Twist property with rotation interval $(0,1)$. For each $\rho\in(0,1)$ there exists a closed, invariant, minimal set ${\cal O}_\rho\subset [0,2\pi)\times(0,\pi)$, which can be injectively projected on $(0,2\pi)$ and such that the induced dynamics preserves the order of $[0,2\pi)\sim S^1$.
This invariant set ${\cal O}_\rho$ can be a periodic orbit,  an Aubry Mather set or a rotational curve , i.e., a continuous closed curve which is the graph of a Lipschitz function.

We will concentrate on billiards having an invariant rotational curve $\gamma=$ graph$(g)$,  $g: [0,2\pi)\mapsto(0,\pi)$ a Lipschitz function.  In general, for a given boundary $\Gamma$, the set of rotation numbers $\rho$ such that ${\cal O}_\rho$ is an invariant rotational curve is nowhere dense in $(0,1)$, but nevertheless invariant rotational curves exist on a large class of classical oval billiards.

For instance, the circular and the elliptical billiards have invariant rotational curves with any rotation number $\rho$. The billiard on a curve of constant width have the horizontal curve $g(\varphi)\equiv \pi/2$ as an invariant rotational curve with rotation number $1/2$. But, apart from special examples like these, it is difficult to find billiards with invariant rotational curves of rational rotation number. In fact,  having an invariant rotational curve with rational rotation number is not a generic property for billiards on ovals. The generic dynamics is to have, for each rational rotation number, a finite number of periodic orbits with this rotation number and at least one of them is hyperbolic, with transverse homoclinic orbits \cite{dia}.

From the other side, having an invariant rotational curve with  irrational rotation number is quite general. If the oval boundary $\Gamma $ is $C^k, k\geq 5$ then  Lazutkin Theorem guarantees a whole family of invariant rotational curves with diophantine $\rho$ near the boundaries of the cylinder $ [0,2\pi)\times(0,\pi)$ (\cite{laz}, \cite{dou}). 
And sufficiently small perturbations of the circular billiard present, in addition to the diophantine invariant rotational curves, given by the KAM theorem, uncountably many invariant rotational curves with liouvillian rotation number.

 For our purposes, we will need the curve $\gamma$, or equivalentely the function $g$, to be $C^2$ or more. This is the case of the circular, the elliptical and the constant width billiards, for instance, where the rotational invariant curves are, in fact, analytic. Also, the curves guaranteed by Lazutkin Theorem are as differentiable as the billiard map, so $C^4$ or more. But, as showed in \cite{arn}, invariant rotational curves can be just a little bit more regular than Lipschitz, and so we have to impose the existence of a $C^2$ rotational invariant curve as a condition to the billiards we will use.

In the next paragraphs we present some bounds that will be necessary to prove our results.

Let us assume that the invariant set ${\cal O}_\rho$ is a closed curve $\gamma$= graph$(g)$. As $B|_\gamma$ preserves the order of $[0,2\pi)\approx S^1$ then either $g(\varphi)\equiv \pi/2$ or there exist  constants $b$ and $B$ such that $0<b\leq g(\varphi)\leq B< \frac{\pi}{2}$ or
$\frac{\pi}{2} <B\leq g(\varphi) \leq b< \pi$.

 The $C^2$ (or more) character of $\gamma$ also implies that a tangent vector $(1,g'(\varphi_0))$ is sent by $DB_{(\varphi_0,\alpha_0)}$ on a tangent vector to $\gamma$ at $(\varphi_1,\alpha_1)$. The preservation of orientation implies that the first coordinate of $DB_{(\varphi_0,\alpha_0)}(1,g'(\varphi_0))$ must be strictly positive.
So  $$L[1+g'(\varphi_0)]-R_0\sin\alpha_0>0 \ \ \hbox{ and } \ \  1+g'(\varphi_0)>\frac{R_0\sin\alpha_0}{L}>0$$
This implies that $g'(\varphi)>-1$ for every $\varphi$. As the billiard is reversible, the graph of $\tilde g(\varphi)=\pi-g(\varphi)$ is also a rotational invariant curve and then $\tilde g'(\varphi)=-g'(\varphi) >-1$
for every $\varphi$. So, for any $C^2$ invariant $\gamma=$graph$(g)$, $-1<g'(\varphi)<1$.

To each invariant rotational curve $\gamma$ is associated a caustic \cite{Tab95},  a curve lying inside the billiard table and
tangent to every segment of the billiard trajectory between two consecutive impacts. If $(\varphi_0,\alpha_0)$ and so
$(\phi_1,\alpha_1)$ belong to $\gamma$, then $\alpha_0=g(\varphi_0)$, $\alpha_1=g(\varphi_1)$ and we have 
\begin{equation}\label{eq:caustica}
\frac{R_0 \sin \alpha_0}{1+g'(\varphi_0)} +  \frac{R_1\sin \alpha_1}{1-g'(\varphi_1)} = L
\end{equation}
where the two terms on the left hand side are strictly positive. The first one measures the distance to the tangency
point (on the caustic) from the initial point $\Gamma(\varphi_0)$  and the second one, from the final point
$\Gamma(\varphi_1)$ (see \cite{cmp}).

\section{Non Elastic Billiards}\label{sec:nonela}

Let $B(\varphi_0,\alpha_0)=(\varphi_1,\alpha_1)$ be a $C^2$ classical billiard map on an oval, with a $C^2$ invariant rotational curve $\gamma_0$, given by the graph of $\alpha=g(\varphi)$.

A compact subset of the phase space $[0,2\pi)\times (0,\pi)$ with non-empty interior and whose boundaries are two distinct rotational curves (not necessarily invariant nor graphs) will be called a compact strip.

Let $I \in \R$ be a closed interval containing $0$. Given $h:I\to I$, a $C^2$ strictly increasing contraction with $h(0)=0$, we can define a non elastic billiard map $P$ on a compact strip $\Sigma$ containing  $\gamma_0$  by
$$P(\varphi_0,\alpha_0) = \left(\varphi_1,\alpha_1-h(\alpha_1-g(\phi_1))\right)$$
with $\Sigma$ chosen such that if $(\varphi,\alpha)\in\Sigma$ then $\alpha-g(\varphi)\in I$.
$P$ is  the composition of a classical billiard followed by a change at the reflection angle, corresponding to a contraction in the vertical fibers of the invariant rotational curve $\gamma_0$.

 Observe that $h(t) \equiv 0$ corresponds to the classical unperturbed billiard and that $h(t) = t$ corresponds to a map that
sends all the points in $\Sigma$ onto the invariant curve (called {\em slap billiard} in \cite{MPS09}).

The derivative of $P$ is given by
\begin{equation}
DP_{(\varphi_0,\alpha_0)} =
\frac{1}{ R_1 \sin \alpha_1}
\left(\begin{array}{cc}
1 & 0   \\
h'_1 g'_1 & 1-h'_1
\end{array}
\right) %\nonumber \\
\left(
\begin{array}{cc}
L -  R_0 \sin \alpha_0  & L   \\
L -  R_0 \sin \alpha_0 -  R_1 \sin \alpha_1 & L -  R_1 \sin \alpha_1
\end{array}
\right) \label{eqn:p}
\end{equation}
where $g'_i=g'(\varphi_i)$ and $h'_i=h'(\alpha_i-g(\varphi_i))$

Our main result is:
\begin{theorem}\label{nonelastic}
Given a classical oval billiard map $B$, with a $C^2$ invariant rotational curve $\gamma_0=\{(\varphi,g(\varphi))\}$, consider a compact strip $\Sigma$ containing $\gamma_0$ and a closed interval
$I \subset \R$, such that $\alpha-g(\varphi)\in I$ if $(\varphi,\alpha)\in\Sigma$.
If $h:I\mapsto \R$ is a $C^2$ function satisfying $h(0)=0$ and $0\leq 1-\underline l<h'(0)<1$
(with $\underline l$ depending only on $\gamma_0$), then there exists a
compact strip $S \subset \Sigma$ such that the non elastic billiard map $P$ defined by $B$, $g$ and $h$ is a
$C^2$-diffeomorphism from $S$ onto $P(S)$. Its limit set $L(P)$ contains $\gamma_0$ and has a dominated splitting. Moreover, the non elastic perturbation does not change the dynamics on  $\gamma_0$.
\end{theorem}

Proof:
The non elastic billiard $P:\Sigma \to P(\Sigma)\subset
[0,2\pi)\times(0,\pi)$ is the composition of the $C^2$ classical
billiard map $B$ with the $C^2$ perturbation of the identity
$(\varphi,\alpha) \mapsto (\varphi,\alpha)-(0, h(\alpha-g(\phi)))$, where
$h$ is a $C^2$ contraction. Then $P$ is a $C^2$ diffeomorphism.

Given
$\delta> 0 $, let $u_{(\varphi,\alpha)} = (1, g'(\varphi) - \delta) $
and $v_{(\varphi,\alpha)} = (1, g'(\varphi) + \delta)$ be two
linearly independent vector fields defining the cone field ${\cal
C}(\varphi,\alpha)$ and the associated quadratic form
$Q_{(\varphi,\alpha)}$ (as in section \ref{sec:cones}).

Using the change of bases matrices, we have
%\begin{eqnarray}[DP_{(\varphi_0,\alpha_0)}]_U =& \label{eqn:dp} \\
%&{\displaystyle \frac{1}{2\delta R_1 \sin \alpha_1}}\left(
%\begin{array}{cc}
%\delta (l_{0} - \delta L) +(1- h'_1) ( \delta l_{1} - l_{01})  &
%\delta (l_{0} + \delta L) -(1- h'_1) ( \delta l_{1} + l_{01}) \\
%\delta (l_{0} - \delta L) - (1- h'_1)( \delta l_{1} - l_{01}) &
%\delta (l_{0} + \delta L) + (1- h'_1) ( \delta l_{1} + l_{01})
%\end{array} \right)\nonumber\end{eqnarray}
{\begin{equation}
\hskip -0.115cm [DP_{(\varphi_0,\alpha_0)}]_U = 
{\displaystyle \frac{1}{2\delta R_1 \sin \alpha_1}}
\left(
\begin{array}{cc}
\delta (l_{0} - \delta L) +(1- h'_1) ( \delta l_{1} - l_{01})  &
\delta (l_{0} + \delta L) -(1- h'_1) ( \delta l_{1} + l_{01}) \\
\delta (l_{0} - \delta L) - (1- h'_1)( \delta l_{1} - l_{01}) &
\delta (l_{0} + \delta L) + (1- h'_1) ( \delta l_{1} + l_{01})
\end{array} 
\right) \label{eqn:dp} 
\end{equation}}
where
\begin{eqnarray*}
l_{0} &=& L(1+g'_0) -  R_0 \sin \alpha_0 \ \ , \ \
l_{1} = L(1-g'_1) -  R_1 \sin \alpha_1 \\
l_{01} &=& L(1+g'_0)(1-g'_1) - (1-g'_1) R_0 \sin \alpha_0 - (1+g'_0)  R_1 \sin \alpha_1 \ .
\end{eqnarray*}

Relation (\ref{eq:caustica}) implies that  for every $(\varphi_0,\alpha_0)$ and $(\varphi_1,\alpha_1)$ on
$\gamma_0$ we have
$$l_{01} =0 \ \ , \ \
l_{0}=\frac{1+g'_0 }{1-g'_1 } R_1 \sin \alpha_1 \ \ , \ \
l_{1}=\frac{1-g'_1 }{1+g'_0 } R_0 \sin \alpha_0 .$$
The billiard boundary $\Gamma$ is an oval and  as it is compact there exist
constants $a$  and  $A$ and a width $D$ such that $0<a\leq
R(\varphi)\leq A$ and $0< L\leq D$. As the invariant curve  is also compact,
for every $(\varphi_0,\alpha_0)$ and
$(\varphi_1,\alpha_1)$
on $\gamma_0$, there exist constants $c$ and $C$ such
that $\displaystyle C\geq l_{0}\geq c>0$ and $\displaystyle C \geq
l_{1} \geq c>0$. So for points on $\gamma_0$ there are constants $0<\underline l\leq 1$ and $0<\overline L$
such that
$$\underline l\leq \frac{l_0}{l_1} \hskip1cm \mbox{and}\hskip1cm 0<\frac{L}{l_1}\leq \overline L
\ .$$
By formula (\ref{eqn:dp}), each entry of the matrix
$[DP_{(\varphi_0,\alpha_0)}]_ U$ for $(\varphi_0,\alpha_0)\in
\gamma_0$ is of the form
$$ \delta
\left(
l_{0} \pm \delta \, L \pm \left(1- h'(0)\right) 
\, l_{1} 
\right)
\geq 
\delta l_1
\left( \frac{l_0}{l_1} - \delta\frac{ L}{l_1} - (1- h'(0)) \right) \geq 
\delta c \left(\underline l -\delta\overline L-(1- h'(0)) \right). $$

Now, if $0\leq 1-\underline l<h'(0)<1$, we can choose $\delta
>0$ such that $\delta c(\underline l -\delta\overline L-(1-
h'(0)))>0$. As $P$ is a $C^2$ diffeomorphism and remembering that
$\gamma_0$ is compact, we can find a strip $S\subset \Sigma$,
containing $\gamma_0$, where $P$ is well defined and all the
entries of $[DP_{(\varphi_0,\alpha_0)}]_ U$ are strictly positive.

Then, by lemma \ref{lema1}, $L(P)\subset S$  has a dominated splitting and contains $\gamma_0$, since $h(0)=0$.
Moreover, $P|_{\gamma_0}=B|_{\gamma_0}$ and the dynamics under $P$, on $\gamma_0$, is the same as under $B$. As $\gamma_0$ and $P$ are $C^2$, for our billiard dynamics, $\gamma_0$ is either a set of periodic points of same period linked by homo/heteroclinic arcs or supports a rational or an irrational rotation.
\hskip .4cm $\blacksquare$

This result will guide us to construct examples of non elastic billiards on ovals with  limit set having a dominated splitting and supporting a rational or an irrational rotation (pieces of type $\mathcal I$ or $\mathcal R$ of Pujals-Sambarino's Theorem). This will be done in the next section.
We will pay attention to the maximal possible size of the strip $S$ and will try to  see if there are other attractors on  $L(P)$  than $\gamma_0$.

\section{Examples}\label{sec:exa}

\subsection{The circle}

The simplest example of a classical billiard with invariant rotational curves is the circular one. This billiard map is linear and is given by $B(\varphi_0,\alpha_0)=(\varphi_0+2\alpha_0,\alpha_0)$. The phase space $[0,2\pi)\times(0,\pi)$ is foliated by invariant horizontal  curves and the dynamics on each one of them is simply a rotation of $2\alpha_0$.

 We pick one invariant  curve $\gamma_0$, defined by $\alpha=g(\varphi)=\overline\beta_0$. At $\gamma_0$,  $g'\equiv 0$, $R_i=R$, the radius of the circle, and $\sin\alpha_i=\sin\overline\beta_0$, implying $l_0=l_1$ and $\underline l=1$.

Fix $I$, a closed interval with $0\in I\subset (-\overline\beta_0,\pi- \overline\beta_0)$ and $h:I\mapsto\R$, any $C^2$ strictly increasing contraction such that $h(0)=0$ and $0< h'(0)<1$. The non elastic billiard $P$ is defined on the strip
$[0,2\pi)\times \{I+\overline\beta_0\}$ and is given by
$P(\varphi_0,\alpha_0)=(\varphi_0+2\alpha_0,\alpha_0-h(\alpha_0-\overline\beta_0))$.

By Theorem~\ref{nonelastic}, there exists a compact strip $S$ such that $P|_S$ is a $C^2$ diffeomorphism and $L(P|_S)$ contains $\gamma_0$ and has a dominated splitting.

Now, we choose $\overline\beta_-$ and $\overline\beta_+$ such that $W=[0,2\pi)\times[\overline\beta_{-},\overline\beta_+]$ is the biggest horizontal straight strip contained in $S$. As each boundary $\gamma_\pm=\{(\varphi, \overline\beta_\pm)\}$ is invariant under $B$, we have that $P(W)\subset W$.
The map $P$ is a horizontal rotation followed by a vertical contraction. Denoting $(\varphi_n,\alpha_n)=P^n(\varphi_0,\alpha_0)$, it is easy to see that $\alpha_n\to \overline\beta_0$, as $n\to+\infty$  and so the horizontal circle $\gamma_0$ is the unique attractor of $P$ in $W$.

 Moreover, the restricted map $P|_{\gamma_0}$ is just a rotation of angle $2\overline\beta_0$.
If $\overline\beta_0/\pi$ is rational, $\gamma_0$ is a normally hyperbolic simple closed curve composed by periodic points of the same period.
If $\overline\beta_0/\pi$ is irrational $\gamma_0$ is a normally hyperbolic closed curve supporting an irrational rotation.

This yields an example of a diffeomorphism, defined on a strip $W$, whose limit set has dominated splitting and is composed by a unique piece of type $\mathcal I$ or $\mathcal R$ of Pujals-Sambarino's Theorem.

Clearly, the size of the strip $W$ depends on the choice of the contraction $h$. Taking for instance $h(x)=\mu x$, with $0<\mu<1$, the non elastic billiard is given by
$P(\varphi_0,\alpha_0)=(\varphi_0+2\alpha_0,\alpha_0-\mu(\alpha_0-\overline\beta_0))$ and the basin of attraction of
$\gamma_0$ contains any straight strip $W=[0,2\pi)\times[\overline\beta_{-},\overline\beta_+]$.

\subsection{The ellipse}

A similar example is given by the classical elliptical billiard. We consider an ellipse $\Gamma$ with eccentricity $e$ and minor axis 1. Its radius of curvature $R$ satisfies $\displaystyle\sqrt{1-e^2}\leq R\leq \frac{1}{1-e^2}$.
The associated  classical billiard map is denoted by $B:[0,2\pi)\times (0,\pi)\mapsto[0,2\pi)\times (0,\pi)$.

This billiard system is integrable: the function
$\displaystyle F(\varphi, \alpha)= \frac{\cos^2\alpha - e^2\cos^2\varphi}{1-e^2\cos^2\varphi}$
is a first integral (see, for instance \cite{berry}) and $[0,2\pi)\times (0,\pi)$ is foliated by the levels of $F=F_0$, with $-\frac{e^2}{1-e^2}<F_0<1$. If $0<F_0<1$, the level set consists of two invariant, analytic and symmetric rotational curves, the lower one contained in $[0,2\pi)\times (0,\frac{\pi}{2})$ and the upper one  in $[0,2\pi)\times (\frac{\pi}{2} ,\pi)$. But, unlike the circular case, $B$ has two elliptic islands of period 2, obstructing the rotational invariant curves to foliate the whole phase-space, as can be seen on Figure 1(left). It also has a hyperbolic 2-periodic orbit, with a saddle connection, corresponding to the level $F_0=0$.

For a fixed $0<F_0<1$ let $\displaystyle\gamma_0$  be the lower invariant rotational  curve in $ F(\phi,\alpha) = F_0 $ (the upper case is analogous). It is the graph of $\alpha=g(\varphi)$ given implicitly by $\cos \alpha = \sqrt{F_0+(1-F_0)e^2\cos^2\phi}.$

We have then that, for any $(\varphi,\alpha)\in \gamma_0$, $\sqrt{(1-F_0)(1-e^2)}\leq\sin\alpha\leq\sqrt{1-F_0}$.
Differentiating twice with respect to $\varphi$ we get  $(1-F_0)e^2\sin 2\phi = g'(\phi) \sin 2\alpha$ and  $2(1-F_0)e^2\cos 2\phi = 2 g'(\phi)  \cos 2\alpha  + g''(\phi) \sin 2\alpha$.
The extremal points  of $g'$ must satisfy $\tan 2 \phi = \tan 2\alpha$ which implies  $\max \{g'(\phi)\} = (1-F_0) e^2 = - \min \{g'(\phi)\}$.
We can then take
\begin{equation}\label{underline}
\underline{l} =
(1-e^2)^2 \left( \frac{1-(1-F_0) e^2}{1+(1-F_0) e^2} \right)^2 \le
\frac{(1+g'_0)^2 R_1 \sin \alpha_1}{(1-g'_1 )^2 R_0 \sin \alpha_0}
= \frac{l_0}{l_1}
\end{equation}
 The associated non elastic billiard map is given by $P(\varphi_0,\alpha_0)=(\varphi_1,\alpha_1-h(\alpha_1-g(\varphi_1))$
where the contraction $h: I\mapsto \R$ is an arbitrary $C^2$ function satisfying $h(0)=0$ and $0<1-\underline l\leq h'(0)$, and $I$ is a closed interval containing $0$. Then, by Theorem~\ref{nonelastic}, there exists a compact strip $S$, containing $\gamma_0$, such that  $P|_S$ is a $C^2$ diffeomorphism,  $L(P|_S)$  has a dominated splitting and contains $\gamma_0$.

Let $F_\pm$ be two constants of motion satisfying  $1> F_- > F_0 > F_+ > 0 $, and $\gamma_\pm$=graph$(g_\pm)$ be the lowest invariant rotational curves associated to $F_\pm$, respectively.
We also suppose that
$W=\{(\varphi,\alpha), g_{-}(\varphi)\leq\alpha\leq g_{+}(\varphi)\}\subset S$ is the biggest strip of this type contained on $S$. As $\gamma_\pm$ are invariant under $B$ we have $P(W)\subset W$.

For $(\varphi_0,\alpha_0)\in W$, we denote $(\varphi_n,\alpha_n)=P^n(\varphi_0,\alpha_0)$. Let us suppose, for instance, that $F(\varphi_0,\alpha_0)>F_0$, the other case being analogous. As $P$ is a translation on a $B$-invariant rotational curve followed by a contraction on the vertical direction toward $\gamma_0$, then $F(\varphi_n,\alpha_n)\to F_0$ and $\gamma_0$ is the $\omega$-limit of any $(\varphi_0,\alpha_0)\in W$. As in the circular case, $\gamma_0$ is the unique attractor of $P|_W$.

Using action-angle variables, Chang and Friedberg \cite{cha} have shown how to decide if the rotation number associated to the level  $F_0=F(\varphi_0,\alpha_0)$ of a given  initial condition $(\varphi_0,\alpha_0)$, is rational or irrational and then determine the dynamics on each level curve.
As $P|_{\gamma_0}=B |_{\gamma_0}$, this allows us to choose $F_0$ in order to have, as the unique attractor, a normally hyperbolic closed curve supporting an irrational rotation or a normally hyperbolic closed curve composed by periodic points of same period.

 Although theoretically promising as a result, depending on the choice of $\gamma_0$ and $h$, the strip $W$ can be very thin. In particular it will never contain points of the elliptical islands.
In trying to go beyond the theoretical predictions and find examples of non elastic elliptical billiards defined on a bigger part of the phase space, we shall remember that weak contractions will never destroy the rotation carried by the linear ellipticity of the 2-periodic orbit. So we can not expect the associated non elastic billiard to be defined on a strip containing the islands. However, it is an interesting question whether this can be achieved by taking a sufficiently strong contraction.

\begin{figure}[h]
\begin{center}
\includegraphics[viewport=0 0 570 570,width=.4\hsize]{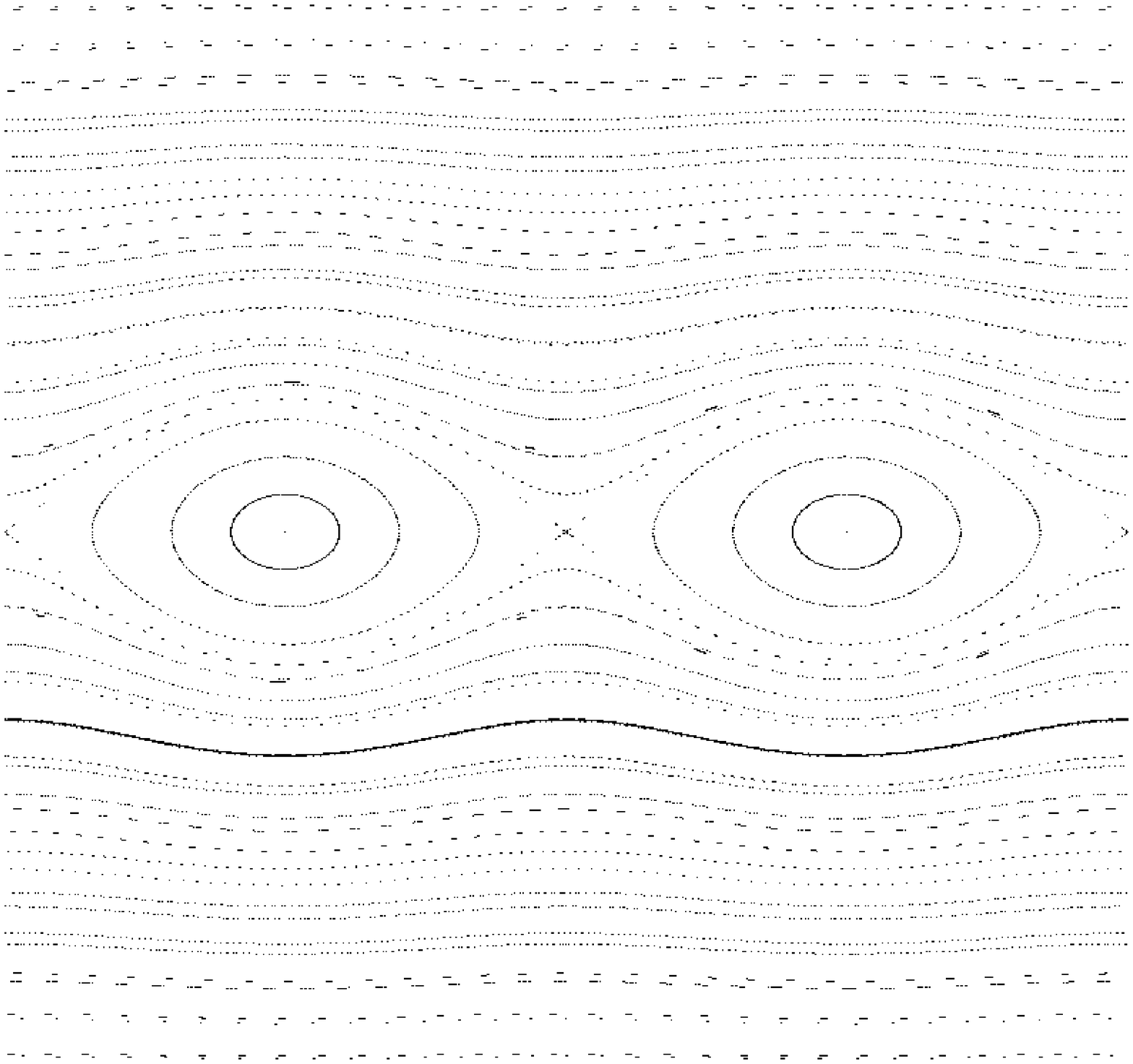}%{elipse-e035-F025}
\hskip 0.5cm
\includegraphics[viewport=0 0 570 570,width=.4\hsize]{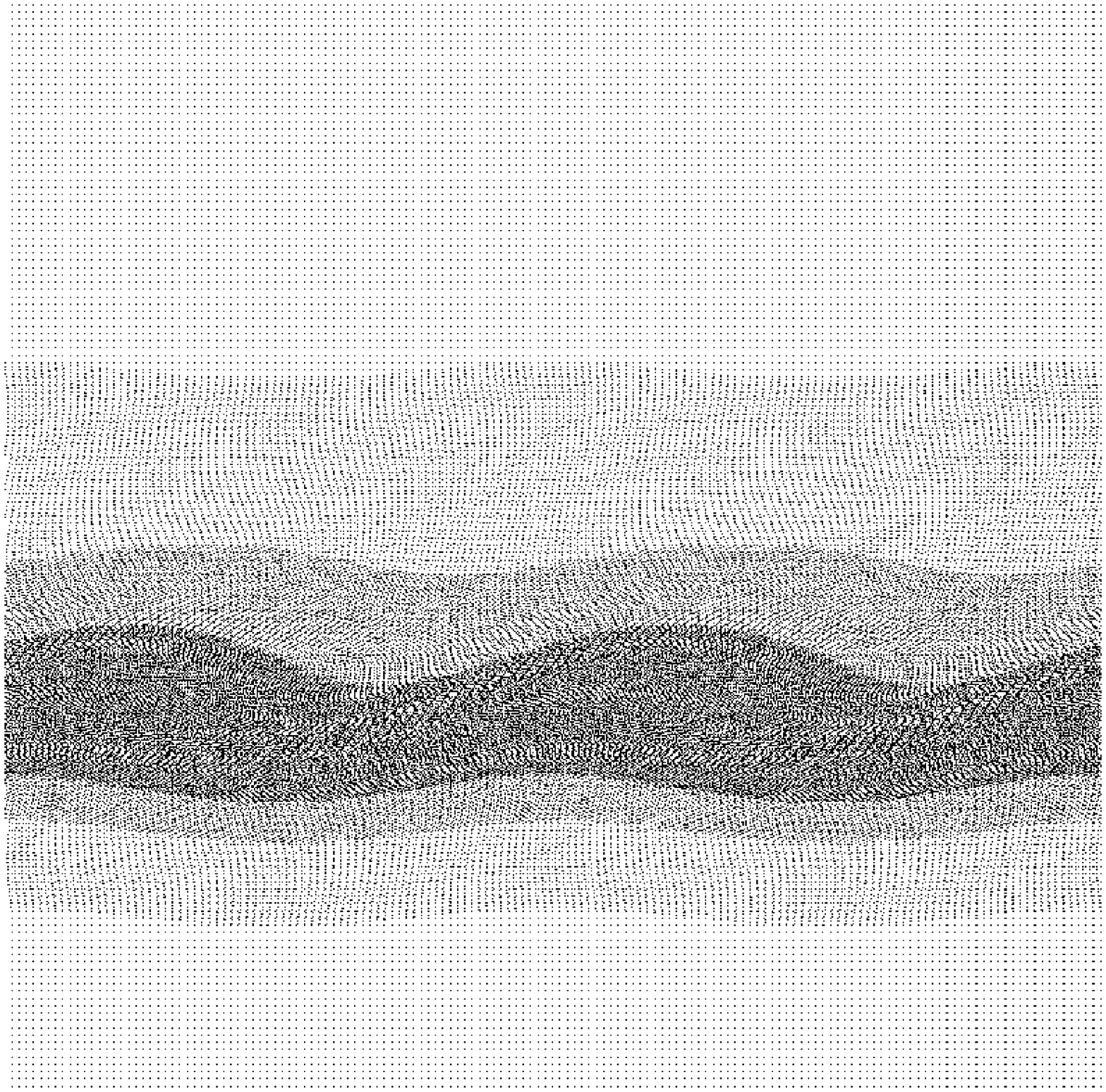}%{elipse-e035-F025-m050-bacia}
\end{center}
\caption{Classical and non elastic elliptical billiards}
\label{fig:elipse1}
\end{figure}

We present some numerical simulations where this can be done. We choose an ellipse with $e= 0.35$ and fix the invariant curve $\gamma_0$ at $F_0=0.25$. Our choices of the eccentricity and $\gamma_0$ are rather arbitrary.
Figure~\ref{fig:elipse1} (left) below diplays the phase space of the classical elliptical billiard. The invariant curve $\gamma_0$ is enhanced. The horizontal axis corresponds to  $\phi \in [0,2\pi)$ and the vertical to $\alpha \in (0,\pi)$, the left bottom corner being the origin.

We consider linear perturbations $h(x) = \mu x$, with $0<\mu<1$.  $\mu =0$ implies that there is no perturbation (classical billiard) and $\mu = 1$ implies that all the points in the phase space land on the invariant curve after one iteration (slap billiard). The non elastic billiard $P(\varphi_0,\alpha_0)=(\varphi_1,\alpha_1-\mu(\alpha_1-g(\phi_1)))$
is a $C^2$ diffeomorphism on any compact strip contained on $[0,2\pi)\times(0,\pi)$. Calculating $\underline l$ by formula (\ref{underline}) we get that   if $\mu>1-\underline l \approx 0.47$
there is a strip $W$ such that $\gamma_0$ is the unique attractor of $P$ on $W$ and has a dominated splitting.
Figure~\ref{fig:elipse1} (right) illustrates the basin of
attraction of the curve $\gamma_0$ for $\mu = 0.5$: black points
correspond to initial conditions which approach $\gamma_0$ under
iteration. This simulation thus indicates that the basin of
attraction of $\gamma_0$ is the whole phase space, ie, that $W$
can be any compact strip contained in $[0,2\pi)\times(0,\pi)$ and
$\gamma_0$ is the unique attractor of $P$. Note that, as we have
mentioned before, either $\gamma_0$ is composed of periodic points
of same period or supports an irrational rotation.

\subsection{Non integrable billiards}

There are no other known $C^2$ ovals such that the classical  billiard map is integrable, other than the circle and the ellipse. But it is well known that there are oval billiards with $C^2$ invariant rotational curves.
Taking a sufficiently differentiable, non integrable, classical oval billiard map $B$ with a $C^2$ invariant rotational curve $\gamma_0$=graph$(g)$ we can, as before, take a contraction $h$ satisfying the hypothesis of Theorem~\ref{nonelastic} and define a non elastic billiard map $P$.
Then there exists a strip $S$ on which $P$ is a $C^2$ diffeomorphism on $S$, $L(P)$ contains $\gamma_0$ and has a dominated splitting. Depending on the rotation number of $\gamma_0$ we can have, as in the previous examples, a normally hyperbolic closed curve supporting an irrational rotation or a normally hyperbolic closed curve composed by periodic points of same period. But we can no longer guarantee that the only attractor of $P$ is $\gamma_0$.

\subsubsection{Invariant straight line}

In order to explore numerically what happens in those more general
cases, we have to pick a concrete example. Although in any example one can prove
the existence of whole families of invariant curves, it is almost
impossible for any fixed curve, to write down the function $g$ for
which it is the graph.
In general this can only be achieved in very specific examples as, for instance, the constant width curves, where $g(\varphi)=\pi/2$. We will deal with a much richer example, given by symmetric perturbations of the circle.

Let $\Gamma_n$ be the oval parameterized by the angle $\varphi$ and
which radius of curvature is of the form $R(\phi)= 1 + a \cos n\phi$, $|a|<1$ and  $n\geq 4$.
Tabachnikov (Section 2.11, \cite{Tab95}) showed that its associated classical billiard map have an invariant
straight line given by $g(\varphi)=\beta_0$ if $\beta_0$ satisfies $n\tan \beta_0 = \tan n\beta_0$. It is not difficult to show that the dynamics on $\gamma_0$ is an irrational rotation \cite{geraldo}.

\begin{figure}[h]
\begin{center}
\includegraphics[viewport=48 0 236 290,width=.15\hsize]{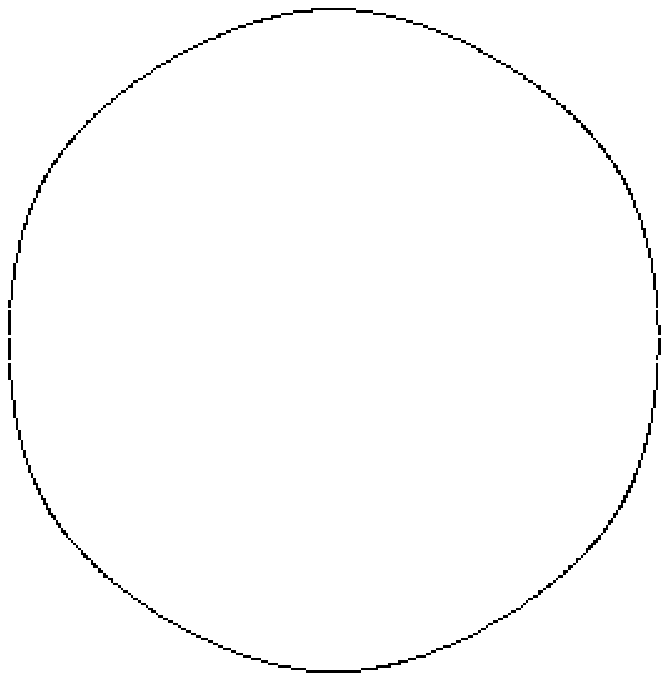}%{sim6-e0010-curva}
\hskip 2cm
\includegraphics[viewport=0 0 570 570,width=.4\hsize]{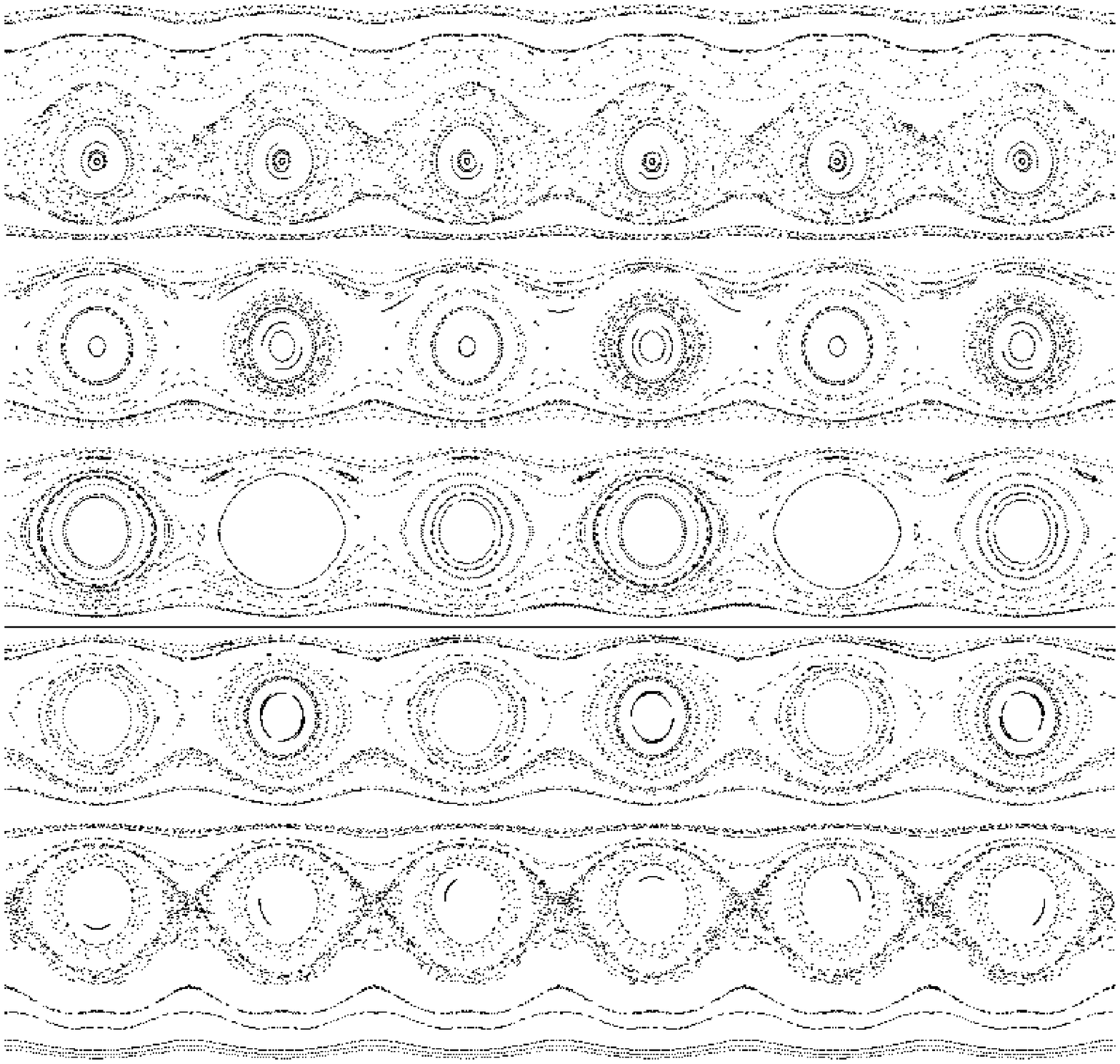}%{sim6-e0010-fase}
\end{center}
\caption{Billiard table and classical phase space for $\Gamma_6$}
\label{fig:mesa}
\end{figure}

For $n=6$, $ \beta_0 = \tan^{-1} \sqrt{7+4\sqrt{21}/3} \approx 0.41\pi$ satisfies this condition and the line $\gamma_0$ given by $\alpha = g(\phi) = \beta_0$  is invariant.
Figure \ref{fig:mesa} displays the billiard table and the corresponding phase space of $\Gamma_6$ defined by $R(\phi) = 1+ 0.01 \cos 6\phi$. As in the previous section, the horizontal axis corresponds to  $\phi \in [0,2\pi)$ and the vertical to $\alpha \in (0,\pi)$, the left bottom corner being the origin. Note the invariant straight line $\gamma_0$.

As in the other numerical examples, we consider linear perturbations $h(x) = \mu x$, $0<\mu<1$.
On $\gamma_0$, we have $l_0 = R_1 \sin \beta_0$ and
$l_1 = R_0 \sin \beta_0$ and we can take $\underline{l} = \frac{\min R(\phi)}{\max R(\phi)} =
\frac{0.99}{1.01}$.
Then, if  $\mu>1-\underline l \approx 0.02$, there exists a strip $S$ on which $P$ is a $C^2$ diffeomorphism, $L(P)$ contains $\gamma_0$ and has a dominated splitting.

\begin{figure}[h]
\begin{center}
\includegraphics[viewport=0 0 570 570,width=.4\hsize]{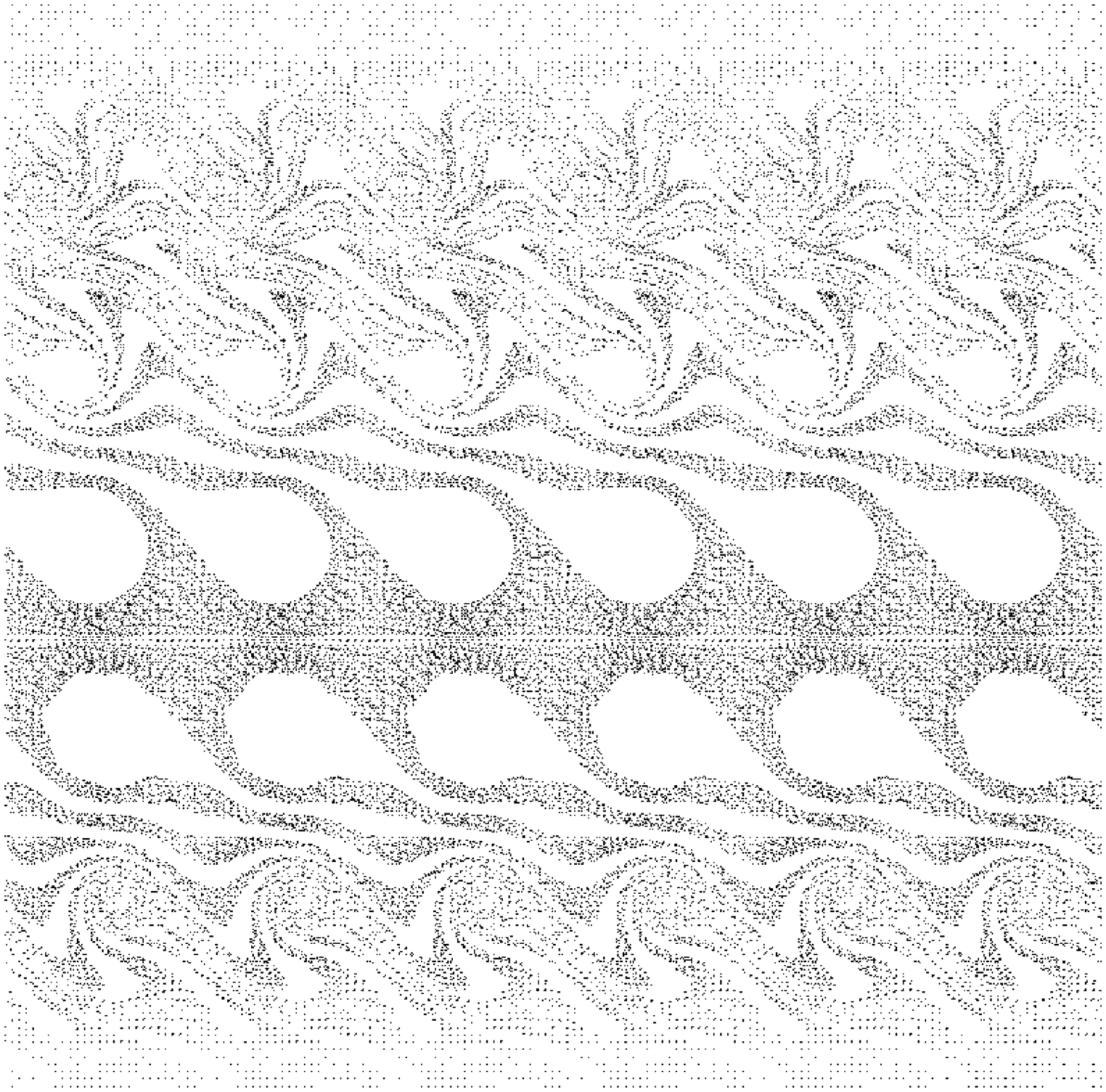}%{sim6-e0010-m010-bacia}
\hskip 0.5cm
\includegraphics[viewport=0 0 570 570,width=.4\hsize]{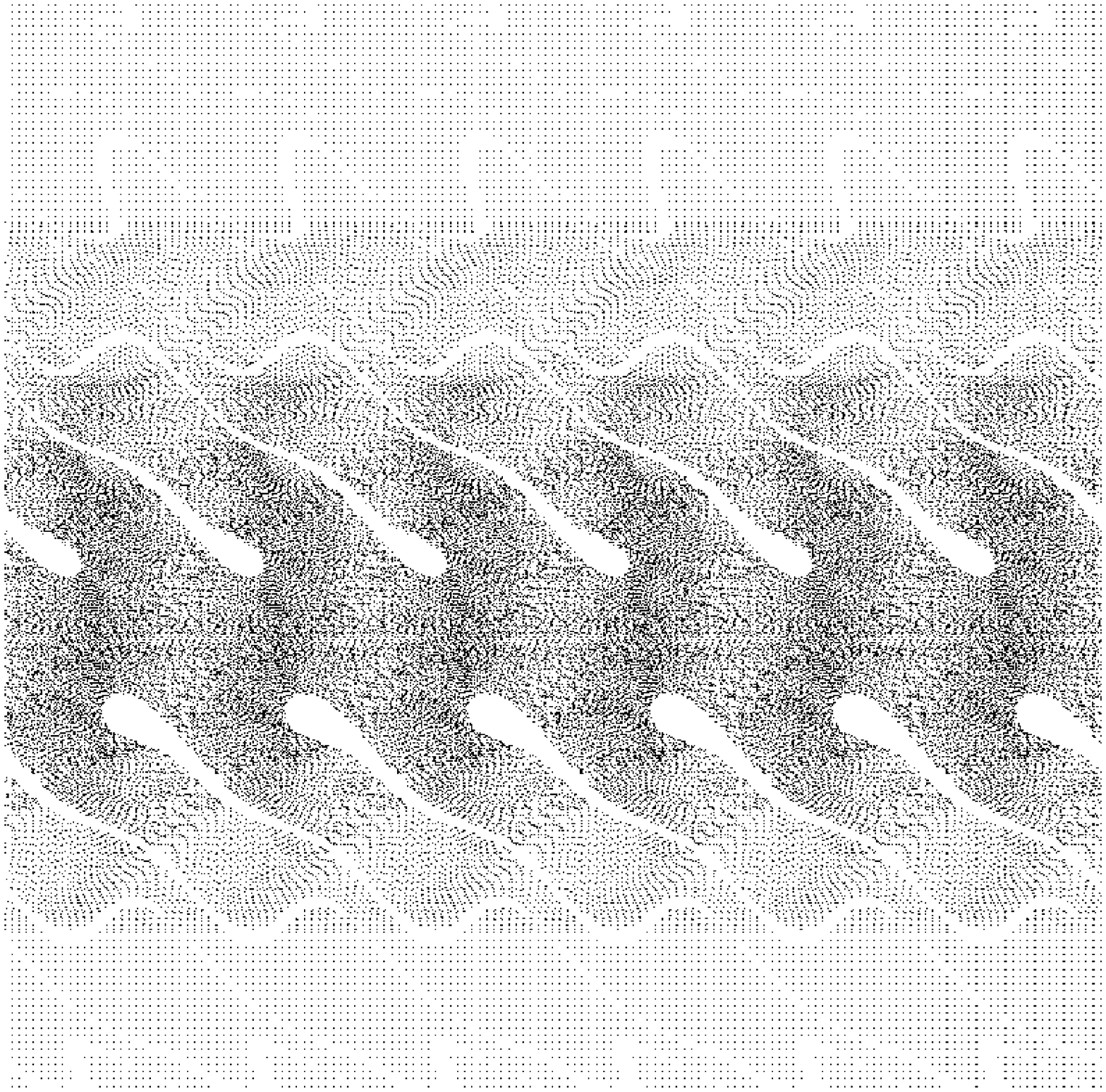}%{sim6-e0010-m035-bacia}

\vskip 0.5cm

\includegraphics[viewport=0 0 570 570,width=.4\hsize]{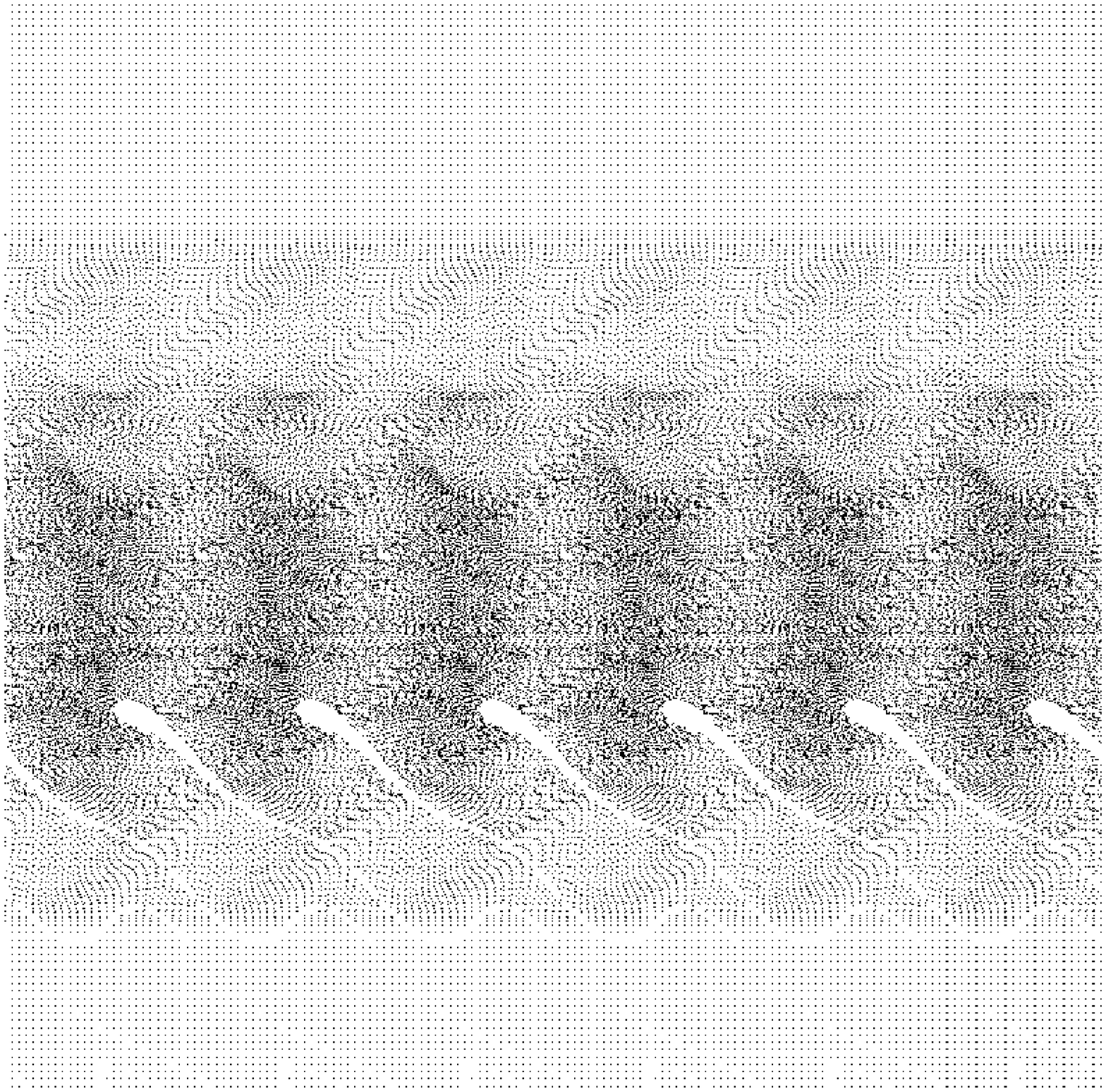}%{sim6-e0010-m037-bacia}
\hskip 0.5cm
\includegraphics[viewport=0 0 570 570,width=.4\hsize]{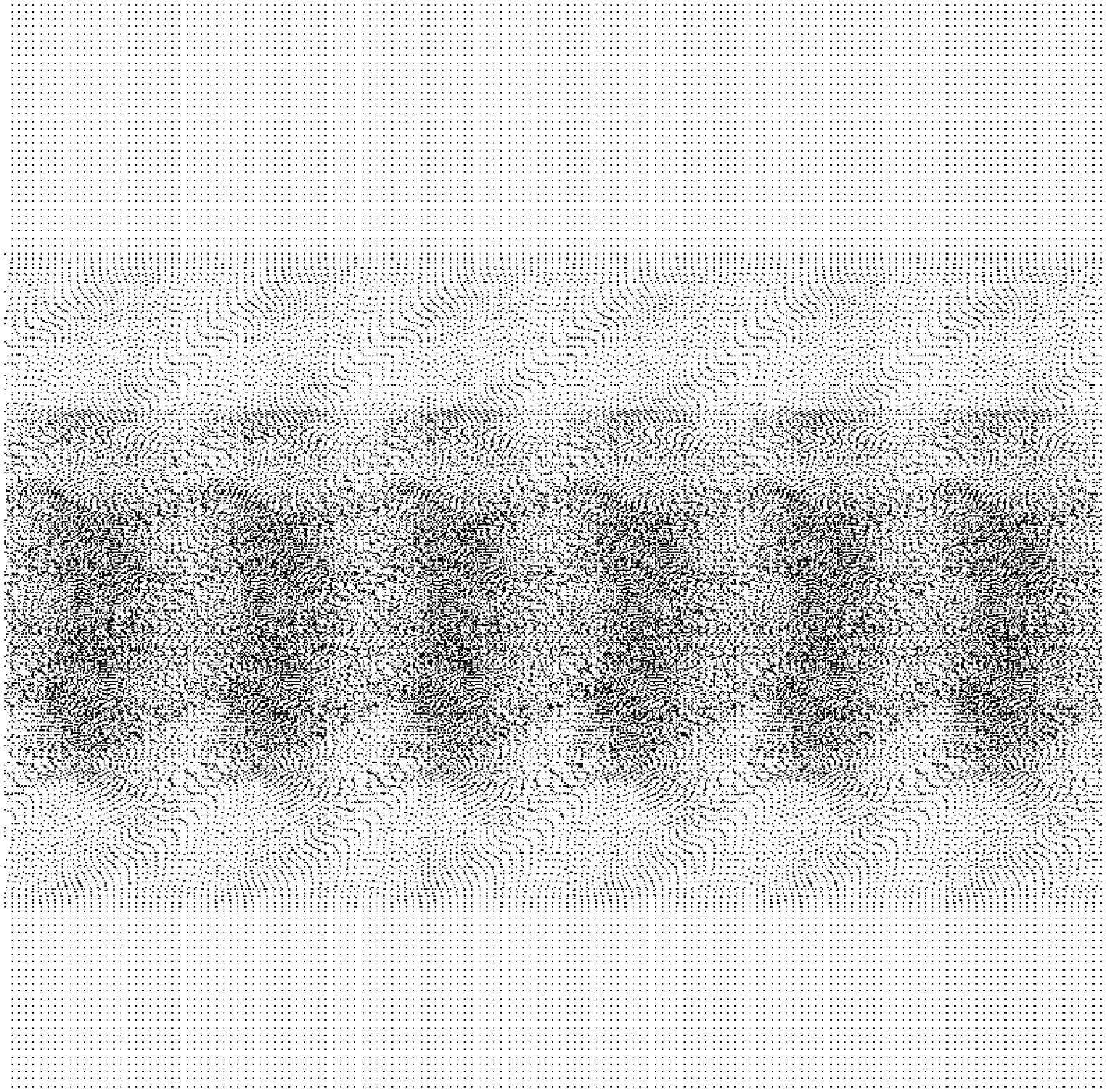}%{sim6-e0010-m040-bacia}
\end{center}
\caption{$\Gamma_6$: the basin of attraction of $\gamma_0$ for $\mu=0.1, 0.35, 0.37$ and $ 0.4$}
\label{fig:4bacias}
\end{figure}

Figure \ref{fig:4bacias}  illustrates the basin of attraction of $\gamma_0$ for $\mu=0.1,\, 0.35,\, 0.37$ and $ 0.4$: black points correspond to initial conditions which approach $\gamma_0$ under iteration. The white tadpoles correspond to points attracted to the 6-periodic orbits, linearly elliptic for the original classical billiard.

We observe that as the contraction factor $\mu$ is increased, the basin of attraction of $\gamma_0$ grows until it eventually occupies the whole phase space.
Thus, for strong contractions, we may have $P$ defined on any compact strip $S$ and having $\gamma_0$ as its unique attractor.

\subsection{An example that is not one}

Any strictly convex classical billiard map is a monotone twist map
with rotation interval $(0,1)$. To each $\rho\in (0,1)$ is
associated a set ${\cal O}_\rho$. If $\rho$ is irrational, ${\cal
O}_\rho$ is either a rotational invariant curve or an Aubry-Mather
set. An Aubry Mather set is a closed, invariant, minimal set,
projecting injectively on a Cantor set of $S^1\sim[0,2\pi)$ and
such that the dynamics preserves the order of $S^1$. It is
contained in a non invariant graph of a continuous piecewise
linear Lipschitz function $\alpha=g(\varphi)$ (see, for instance \cite{kat}, Section 13.2).

Let us take a classical billiard map $B$ with two rotational invariant curves $\gamma_{-}$ and $\gamma_+$ and an Aubry-Mather set ${\cal A}$ contained in the strip bounded by $\gamma_{-}$ and
$\gamma_+$ and a perturbation $h$ on this strip.
The non elastic billiard map can be then defined as before, as $P(\varphi_0,\alpha_0) =(\varphi_1, \alpha_1 -h(\alpha_1 - g(\varphi_1))$, where $(\varphi_1, \alpha_1)=B(\varphi_0, \alpha_0)$.

As $h$ is a contraction, the limit set of $P$ will contain the Aubry-Mather set ${\cal A}$. However, we must remark that $P$ is not $C^2$ because $g$ is not even a differentiable function. Then, we can not apply Pujal-Sambarino's theorem.

\noindent {\bf Acknowledgment.} 
{We thank M\'ario Jorge Dias Carneiro for his precious ideas. We thank the
Brazilian agencies FAPEMIG and CNPq for financial support.
}

\end{document}